\documentclass[12pt,a4paper]{elsart}
\usepackage{amsmath,amssymb}
\usepackage{natbib}
\usepackage{color}
\journal{jspi}

\newenvironment{proof}{\rm \trivlist  \item[\hskip \labelsep{\sc Proof.}]}
{\qed \endtrivlist}

\newtheorem{definition}{Definition}
\newtheorem{corollary}{Corollary}
\newtheorem{proposition}{Proposition}

\usepackage{amsmath,amssymb}
\usepackage[dvips]{epsfig}
\newcommand{\design}{{\mathcal D}}
\newcommand{\fraction}{\mathcal F}

\newcommand{\jei}{\text{i }}
\def\cocoa{{\hbox{\rm C\kern-.13em o\kern-.07em C\kern-.13em o\kern-.15em A}}}
\begin{document}
\begin{frontmatter}

\title{Indicator function and complex coding for mixed fractional factorial designs\thanksref{grant}}%

\author{Giovanni Pistone}
\address{Department of Mathematics - Politecnico di Torino}
\ead{giovanni.pistone@polito.it}
\author{Maria-Piera Rogantin\corauthref{cor}}
\corauth[cor]{Corresponding author.}
\address{Department of Mathematics - Universit\`a di Genova}
\ead{rogantin@dima.unige.it}
\thanks[grant]{Partially supported by the Italian  PRIN03 grant coordinated by G. Consonni}

\maketitle

\begin{abstract}
In a general fractional factorial design, the $n$-levels of a factor are coded by the $n$-th roots of the unity. This device allows a full generalization to mixed-level designs of the theory of the polynomial indicator function which has already been introduced for two level designs by Fontana and the Authors (2000). the properties of orthogonal arrays and regular fractions are discussed.
\end{abstract}


\end{frontmatter}
\section{Introduction}
Algebraic and geometric methods are widely used in the theory of the design of experiments. A variety of these methods
exist: real linear algebra, $\mathbb Z_p$ arithmetic, Galois Fields GF$(p^s)$ arithmetic, where $p$ is a prime number
as in \cite{bose:47}. See, e.g., \cite{raktoe|hedayat|federer:81} and the more recent books by \cite{dey|mukerjee:1999}
and \cite{wu|hamada:2000}.

Complex coding of levels has been used by many authors in various contexts, see e.g. \cite{bailey:1982},
\cite{kobilinsky|monod:91}, \cite{edmondson:94}, \cite{kobilinsky|monod:95}, \cite{collombier:96} and
\cite{xu|wu:2001}.

The use of a new background, called Commutative Algebra or Polynomial Ring Algebra, was first advocated by
\cite{pistone|wynn:96} and later discussed in detail in \cite{pistone|riccomagno|wynn:2001}. Other relevant general
references are \cite{robbiano:98}, \cite{robbiano|rogantin:98} and \cite{galetto|pistone|rogantin:2003}.

In the present paper, mixed-level (or asymmetric) designs with replicates are considered and  the approach to the
two-level designs discussed in \cite{fontana|pistone|rogantin:97} and \cite{fontana|pistone|rogantin:2000} is
generalized. In the latter, the fractional factorial design was encoded in its indicator function with respect to the
full factorial design. In \cite{tang|deng:99}, entities related to coefficients of the polynomial indicator function were
independently introduced into the construction of a generalized word length pattern. The coefficients themselves were
called $J$-characteristics in \cite{tang:2001}, where it was shown that a two-level fractional design is uniquely
determined by its $J$-characteristics. The representation of a fraction by its indicator polynomial function was
generalized to designs with replicates in \cite{ye:2003} and extended to non two-level factors using orthogonal
polynomials with an integer coding of levels in \cite{cheng|ye:2004}.

Sections 2 and 3 are a self-contained introduction of the indicator function representation of a factorial design using
complex coding. The main results are in Sections 4 to 6.  The properties of the indicator polynomial are discussed in
Section \ref{sec:fraction}. If the factor levels are coded with the $n$-th roots of the unity, the coefficients of the
indicator polynomial are related  to many interesting properties of the fraction in a simple way: orthogonality among the
factors and interactions, projectivity, aberration and regularity. Combinatorial orthogonality vs. geometrical
orthogonality is discussed in Section \ref{sec:project}. A type of generalized regular fraction is defined and discussed
in Section \ref{sec:regular}. The usual definition, where the number of levels is prime for all factors is extended  to
asymmetric design with any number of levels. With such a definition, all the monomial terms of any order are either
orthogonal or totally aliased. However, our framework does not include the GF($p^s$) case. Some examples are shown in
Section~\ref{sec:ex}.

A first partial draft of the present paper was presented in the GROSTAT V 2003 Workshop. Some of the results of
Proposition \ref{th:reg} have been obtained independently by \cite{ye:2004}.

\section{Coding of factor levels}

Let $m$ be the number of factors of a design. We denote the factors by $A_j$, $j=1,\dots,m$, and  the number of levels of
the factor $A_j$ by $n_j$. We consider only \emph{qualitative factors}.

We denote the full factorial design by $\design$, $\design = A_1 \times \cdots \times A_m$, and the space of all real
responses defined on $\design$ by $\mathcal{R}(\design)$.

In some cases, it is of interest to \emph{code} qualitative factors with numbers, especially when the levels are ordered.
Classical examples of numerical coding with rational numbers $a_{ij} \in \mathbb Q$ are: (1) $a_{ij}=i$, or (2)
$a_{ij}=i-1$, or (3) $a_{ij}=(2i-n_j-1)/2$ for odd $n_j$  and $a_{ij}=2i-n_j-1$ for even $n_j$, see \cite[Tab.
4.1]{raktoe|hedayat|federer:81}. The second case, where the coding takes value in the additive group $\mathbb Z_{n_j}$,
i.e. integers mod $n_j$, is of special importance.  We can define the important notion of \emph{regular fraction} in such
a coding. The third  coding is the result of the orthogonalization of the linear term in the second coding with respect
to the constant term. The coding $-1, +1$ for two-level factors has a further property in that the values $-1, +1$ form a
multiplicative group. This property was widely used in \cite{fontana|pistone|rogantin:2000}, \cite{ye:2003},
\cite{tang|deng:99} and \cite{tang:2001}.

In the present paper, an approach  is taked to parallel our theory for two-level factors with coding $-1, +1$. The $n$
levels of a factor are coded by the complex solutions of the equation $\zeta^n=1$:
\begin{equation} \label{omega-k}
\omega_k=\exp\left(\jei \frac {2 \pi}{n} k \right) \quad , \qquad \ k=0,\ldots, n-1 \ .
\end{equation}
We denote  such a factor with $n$ levels by $\Omega_n$, $\Omega_n=\left\{\omega_0,\ldots,\omega_{n-1}\right\}$. With such
a coding, a complex orthonormal basis of the responses on the full factorial design is formed by all the monomials.

For a basic reference  to the algebra of the complex field $\mathbb C$ and of the $n$-th complex roots of the unity
references can be made to \cite{lang:65}; some useful points are collected in Section \ref{pr:roots} below.

As $\alpha=\beta \mod n$ implies $\omega_k^\alpha = \omega_k^\beta$, it is useful to introduce the residue class ring
$\mathbb Z_n$ and the notation $[k]_n$ for the residue of $k \mod n$. For integer $\alpha$, we obtain $(\omega_k)^\alpha
= \omega_{[\alpha k]_n}$. The mapping
\begin{equation}\label{map}
\mathbb Z_n  \longleftrightarrow  \Omega_n \subset \mathbb C  \qquad \text{ with } \qquad
 k  \longleftrightarrow \omega_k
\end{equation}
is a group isomorphism of the additive group of $\mathbb Z_n$ on the multiplicative group  $\Omega_n \subset \mathbb C$.
In other words,
\begin{equation*}
\omega_h\omega_k=\omega_{[h+k]_n} \ .
\end{equation*}
We drop the sub-$n$ notation when there is no ambiguity.

We denote by:
 \begin{itemize}
 \item   ${\#\design}$: the number of points of the full factorial design, ${\#\design}=\prod_{j=1}^m n_j$.
 \item   $L$:  the full factorial design with integer coding
$\{0,\ldots,n_j-1\}$, $j=1,\dots,m$, and $\design$  the full factorial design with complex coding:
\begin{equation*}
L  = \mathbb Z_{n_1} \times \cdots  \times \mathbb Z_{n_m} \qquad  \textrm{and} \qquad \design =
\design_1 \times \cdots \design_j \cdots
\times \design_m \ \textrm{ with } \design_j=\Omega_{n_j}\end{equation*}
 According to  map (\ref{map}), $L$ is both the integer coded design and the  exponent set of the complex coded design;
 \item  $\alpha$, $\beta$, \ldots: the elements of $L$:
\begin{equation*}
L = \left\{ \alpha= (\alpha_1,\ldots,\alpha_m) : \alpha_j = 0,\ldots,n_j-1, j=1,\ldots,m\right\} \ ;
\end{equation*}
that is, $\alpha$ is both a treatment combination in the integer coding and a multi-exponent of an interaction term;
 \item  $[\alpha-\beta]$:
the $m$-tuple $\left(\left[\alpha_1-\beta_1 \right]_{n_1}, \ldots, \left[\alpha_j-\beta_j \right]_{n_j}, \ldots,
\left[\alpha_m - \beta_m\right]_{n_m} \right)$; the computation of the $j$-th element is in the ring $\mathbb Z_{n_j}$.
\end{itemize}

\section{Responses on the design} \label{sec:bas-des}

The responses on the design and the linear models are discussed in this section. According to the generalization of the
algebraic approach by \cite{fontana|pistone|rogantin:2000}, the design $\design$ is identified as the zero-set of the
system of polynomial equations
\begin{equation*}
\zeta^{n_j}_j-1=0 \quad , \qquad  j=1,\ldots,m \ .
\end{equation*}
A \emph{complex} response $f$ on the design $\design$ is a $\mathbb C$-valued function defined on $\design$. This
response can be considered as the restriction to $\design$ of a complex polynomial.

We denote by:
\begin{itemize}
\item $X_i$; the $i$-th component function, which maps a  point to its $i$-th component:
\begin{equation*}
X_i : \quad \design \ni (\zeta_1,\ldots,\zeta_m)\ \longmapsto \
\zeta_i \ .
\end{equation*}
The function $X_i$ is called \emph{simple term} or, by abuse of terminology, \emph{factor}.
 \item $X^\alpha$, with
$\alpha\in L$: the \emph{interaction term} $X_1^{\alpha_1} \cdots X_m^{\alpha_m}$, i.e. the
function
\begin{equation*}
X^\alpha : \quad \design \ni (\zeta_1,\ldots,\zeta_m)\ \mapsto \
\zeta_1^{\alpha_1}\cdots \zeta_m^{\alpha_m} \quad , \qquad \alpha \in
L \ .
\end{equation*}
The function $X^\alpha$ is a special response that we call \emph{monomial response} or
\emph{interaction term},  in analogy with current terminology.
\end{itemize}

In the following, we shall use the word \emph{term} to indicate either a simple term or an interaction term.

We say term $X^\alpha$ has \emph{order} (or order of interaction) $k$ if $k$ factors are involved,
i.e. if the $m$-tuple $\alpha$ has $k$ non-null entries.

If $f$ is a response defined on $\design$ then its \emph{mean value} on $\design$, denoted by $E_{\design}(f) $, is:
\begin{equation*}
E_{\design}(f)  = \frac 1 {\#\design} \sum_{\zeta \in \design} f(\zeta) \ .
\end{equation*}
We say that a response $f$ is \emph{centered} if $E_\design(f) = 0$. Two responses $f$ and $g$ are \emph{orthogonal on
$\design$} if $E_\design(f \ \overline g) = 0$.

it should be noticed that the set of all the responses is a complex Hilbert space with the Hermitian product $ f\cdot g =
E_\design(f \ \overline g) $.

Two basic properties connect the algebra to the Hilbert structure, namely
\begin{enumerate}
 \item $X^\alpha\overline{X^\beta}=X^{[\alpha-\beta]}$;
 \item $E_{\design}(X^0)=1$, and $E_{\design}(X^\alpha)=0$ for $\alpha
 \neq 0$, see Section \ref{pr:roots} Item (\ref{it:s-poly}).
\end{enumerate}

The set of functions $\left\{X^\alpha \ , \ \alpha \in L \right \}$ is an orthonormal basis of the \emph{complex
responses} on  design $\design$. From  properties (1) and (2) above it follows that:
\begin{equation*}E_{\design}(X^\alpha\overline{X^\beta}) = E_{\design}(X^{[\alpha-\beta]}) =
\begin{cases}
1 & \text{if }  \alpha=\beta \\
0 & \text{if }  \alpha \neq \beta
\end{cases}\end{equation*}
Moreover, $\#L=\#\design$.

Each response $f$ can therefore be represented as a unique $\mathbb C$-linear combination of constant, simple and
interaction terms:
\begin{equation} \label{eq:f-repr}
f = \sum_{\alpha \in L} \theta_\alpha \ X^\alpha, \quad \theta_\alpha
\in \mathbb C
 \end{equation}
where the coefficients are uniquely defined by: $ \theta_\alpha=E_\design \left(f
\overline{X^\alpha}\right) $.
 In fact,
 \begin{equation*}
 \sum_{\zeta \in \design} f(\zeta)
\overline{X^\alpha}(\zeta)= \sum_{\zeta \in \design} \sum_{\beta \in L} \theta_\beta X^\beta \overline{X^\alpha}(\zeta)= \sum_{\beta \in L}
\theta_\beta \sum_{\zeta \in \design}X^\beta(\zeta) \overline{X^\alpha}(\zeta)= {\#\design} \ \theta_\alpha \ .
 \end{equation*}
We can observe that a function is centered on $\design$ if, and only if, $\theta_0=0$.

As $\overline{\theta_\alpha}=E_\design \left(\overline f {X^\alpha}\right)$, the conjugate of response $f$ has the
representation:
\begin{equation*}
\overline{f(\zeta)}=  \sum_{\alpha \in L} \overline{\theta_\alpha} \  \overline{X^\alpha}(\zeta)= \sum_{\alpha \in L}
\overline{\theta_{[-\alpha]}} {X^\alpha}(\zeta) \ .
\end{equation*}
A response $f$ is  \emph{real valued} if, and only if, $ \overline{\theta_\alpha} = \theta_{[-\alpha]}$  for all $\alpha
\in L$.

We suggest the use of the roots of the unity because of the mathematical convenience we are going to show. In most of the
applications, we are interested in \emph{real valued} responses, e.g. measurements, on the design points. Both the real
vector space $\mathcal{R}(\design)$ and the complex vector space $\mathcal{C}(\design)$ of the responses on the design
$\design$ have a real basis, see \cite[Prop. 3.1]{kobilinsky:90} and \cite{pistone|rogantin:2005-30}, where a special
real basis that is common to both spaces is computed. The existence of a real basis implies the existence of real linear
models even though the levels are complex.

\section{Fractions} \label{sec:fraction}
A fraction $\fraction$ is a subset of the design, $\fraction \subseteq \design$. We can algebraically describe a fraction
in two ways, namely using generating equations or the indicator polynomial function.
\subsection{Generating equations} \label{sec:gen-eq}
All fractions can be obtained by adding further polynomial equations, called \emph{generating equations}, to the design
equations $X_j^{n_j}-1 = 0$, for $j=1,\ldots,m$,  in order to restrict the number of solutions.

For example, let us consider a classical $3^{4-2}_{\text{III}}$ regular fraction, see \cite[Table 5A.1]{wu|hamada:2000},
coded with complex numbers  according to the map in Equation (\ref{map}). This fraction  is defined by $ X^{3}_j-1=0$ for
$j=1,\ldots,4 $, together with the generating equations $X_1 X_2 X_3^2 = 1$ and $ X_1 X_2^2 X_4 = 1$. Such a
representation of the fraction is classically termed ``multiplicative'' notation. In our approach, it is not a question
of notation or formalism, but rather the equations are actually defined on the complex field $\mathbb C$. As the recoding
is a homomorphism from the additive group $\mathbb Z_3$ to the multiplicative group of $\mathbb C$, then the additive
generating equations in $\mathbb Z_3$ (of the form $A +B+2C= 0 \mod 3$ and $A+2B+D=0 \mod 3$) are mapped to the
multiplicative equations in $\mathbb C$. In this case, the generating equations are binomial, i.e. polynomial with two
terms.

In the following, we consider general subsets of the full factorial design and, as a consequence, no special form of the
generating equations is assumed.
\subsection{Responses defined on the fraction, indicator and counting functions}

The indicator polynomial  was first introduced in \cite{fontana|pistone|rogantin:97} to describe a fraction. In the
two-level case, \cite{ye:2003} suggested  generalizing the idea of indicator function to fractions with replicates.
However, the single replicate case has special features, mainly because, in such a case, the equivalent description with
generating equations is available. For coherence with  general mathematical terminology, we have maintained the indicator
name, and introduced the new name, that is,  \emph{counting function} for the replicate case. The design with replicates
associated to a counting function  can be considered a multi-subset $\fraction$ of the design $\design$, or an array with
repeated rows. In the following, we also use the name ``fraction''  in this extended sense.

\begin{definition}[Indicator function and counting function]
\label{de:counting}

 The \emph{counting function} $R$ of a fraction $\fraction$ is a response defined on $\design$ so that for
each $\zeta \in \design$, $R(\zeta)$ equals the number of appearances of $\zeta$ in the fraction.

A 0-1 valued counting function is called \emph{indicator function} $F$ of a single replicate fraction $\fraction$.

We denote  the coefficients of the representation of $R$ on $\design$ using  the monomial basis by $b_\alpha$:
\begin{equation*}
R(\zeta)=\sum_{\alpha \in L} b_\alpha \ X^\alpha(\zeta) \qquad \zeta\in\design \  .
\end{equation*}
\end{definition}

A polynomial function $R$ is a counting function of some fraction $\fraction$ with replicates up to $r$ if, and only if,
$R(R-1) \cdots (R-r) = 0$ on $\design$. In particular a function $F$ is an indicator function  if, and only if, $F^2 - F
= 0$  on $\design$.

If $F$ is the indicator function of the fraction $\fraction$, $F-1=0$ is a set of generating equations of the same
fraction.

As  the counting function is real valued, we obtain $\overline{b_\alpha} = b_{[-\alpha]}$.

If $f$ is a response on $\design$ then its \emph{mean value on  $\fraction$}, denoted by $E_{\fraction}(f)$, is:
\begin{equation*}
E_{\fraction}(f)= \frac 1 {\# \fraction} \sum_{\zeta\in\fraction} f(\zeta) = \frac {\#\design} {\#
\fraction} \ E_{\design} (R \ f)
\end{equation*}
where $\# \fraction$ is the total number of treatment  combinations of the fraction, $\# \fraction=\sum_{\zeta \in
\design} R(\zeta)$.

\begin{proposition} \label{pr:bc-alpha}
\begin{enumerate}
 \item
The coefficients $b_\alpha$ of the the counting function of a fraction $\fraction$ are:
\begin{equation*} \label{bc-alpha}
b_\alpha= \frac 1 {\#\design} \sum_{\zeta \in \fraction} \overline{X^\alpha(\zeta)} \ ;
\end{equation*}
in particular, $b_0$ is the ratio between the number of points of the fraction and those of the design.
 \item
In a single replicate fraction, the coefficients $b_\alpha$ of the indicator function are related according to:
\begin{equation*} \label{b-alpha2}
b_\alpha = \sum_{\beta \in L} b_\beta \ b_{[\alpha - \beta]} \ .
\end{equation*}
 \item
If $\fraction$ and $\fraction'$ are complementary fractions without replications and $b_\alpha$ and $b_\alpha'$ are the
coefficients of the respective indicator functions,  $b_0=1-b_0'$ and $b_\alpha = - b_\alpha'$.
\end{enumerate}
\end{proposition}
\begin{proof}
Item (1) follows from :
$$
 \sum_{\zeta \in \fraction} \overline{X^\alpha(\zeta)}  \sum_{\zeta \in \design} R \
  \overline{X^\alpha(\zeta)}= \sum_{\zeta \in \design}
 \sum_{\beta \in L} b_\beta X^\beta(\zeta)   \overline{X^\alpha(\zeta)}
   \sum_{\zeta \in \design}  b_\alpha =
 {\#\design} \ b_\alpha \ .
 $$
Item (2) follows from relation $F=F^2$. In fact:
\begin{displaymath}
\begin{split}
 \sum_\alpha b_\alpha X^{\alpha} & = \sum_\beta b_\beta X^{\beta} \sum_\gamma b_\gamma X^{\gamma }
 = \sum_{\beta, \gamma} b_\beta b_\gamma X^{[\beta+\gamma]} = \\
 & = \sum_\alpha \sum_{[\beta+\gamma]=\alpha} b_\beta b_\gamma X^{\alpha}= \sum_\alpha \sum_\beta  b_\beta \ b_{[\alpha - \beta]} X^{\alpha}\ .
\end{split}
\end{displaymath}
Item (3)  follows from $F'=1-F$.
\end{proof}

\subsection{Orthogonal responses on a fraction} \label{sec:orth}
In this section, we discuss the general case of fractions $\fraction$ with or without replicates. As in the full design
case, we say that a response $f$ is \emph{centered} on a fraction $\fraction$ if $E_\fraction(f) =E_\design(R \ f)= 0$
and we say that two responses $f$ and $g$ are \emph{orthogonal} on $\fraction$ if $E_\fraction(f \ \overline
g)=E_\design(R \ f \ \overline g) = 0$, i.e. the response $f \ \overline g$ is centered.

It should be noticed that the term ``orthogonal'' refers to vector orthogonality with respect to a given Hermitian
product. The standard practise  in orthogonal array literature, however, is to define an array as orthogonal when all the
level combinations appear equally often in relevant subsets of columns, e.g. \cite[Def.
1.1]{hedayat|sloane|stufken:1999}. Vector orthogonality is affected by the coding of the levels, while the definition of
orthogonal array is purely combinatorial. A characterization of orthogonal arrays can be based on vector orthogonality of
special responses.  This section and the next one are devoted to discussing how the choice of complex coding makes such a
characterization as straightforward as in the classical two-level case with coding -1,+1 .

\begin{proposition}\label{pr:ort-bc-alpha} \
Let $R = \sum_{\alpha \in L} b_\alpha X^\alpha$ be the counting function of a fraction~$\fraction$.
\begin{enumerate}
 \item The term $X^\alpha$ is centered on $\fraction$ if, and only if, $b_\alpha=b_{[-\alpha]}=0$.
 \item The terms $X^\alpha$ and $X^\beta$
are orthogonal on $\fraction$ if, and only if, $b_{[\alpha-\beta]}=0$;
 \item If $X^\alpha$ is centered then, for each $\beta$ and $\gamma$ such that $\alpha=[\beta-\gamma]$ or
 $\alpha=[\gamma-\beta]$, $X^\beta$ is orthogonal to $X^\gamma$.
 \item
\label{ort-4} A fraction $\fraction$ is self-conjugate, that is, $R(\zeta) = R(\overline \zeta)$ for any $\zeta \in
\design$, if, and only if, the coefficients $b_\alpha$ are real for all $\alpha \in L$.
\end{enumerate}
\end{proposition}
\begin{proof} \
The first three Items follow easily from Proposition \ref{pr:bc-alpha}.

For the  Item (4), we obtain:
\begin{equation*}
\begin{split}
R(\zeta)= & \sum_{\alpha \in L} b_\alpha X^\alpha(\zeta)=\sum_{\alpha \in L}b_{[-\alpha]} X^{[-\alpha]}(\zeta)
=\sum_{\alpha \in L} \overline{b_\alpha} X^{[-\alpha]}(\zeta)\\
R(\overline\zeta)= & \sum_{\alpha \in L} b_\alpha X^\alpha(\overline\zeta)=\sum_{\alpha \in L}b_\alpha
X^{[-\alpha]}(\zeta)  \ .
\end{split}
\end{equation*}
Therefore  $R(\zeta)=R(\overline{\zeta})$ if, and only if, $b_\alpha= \overline{b_{\alpha}}$. It should be noticed that
the same applies to all real valued responses.
\end{proof}

Interest in self-conjugate fractions concerns the existence of a real valued linear basis of the response space, as
explained in \cite[Prop. 3.1]{kobilinsky:90}. It follows that it is possible to fit  a real linear model on such a
fraction, even though  the levels have complex coding.

An important property of the centered responses follows from the structure of the roots of the unity as a cyclical group.
This connects the combinatorial properties to the coefficients $b_\alpha$'s through the following two basic properties
which hold true for the full design $\design$.
\begin{enumerate}
 \item [P-1]
Let $X_i$ be a simple term with level set $\Omega_n$.  Let us define $s=n/\text{gcd}(r,n)$ and let $\Omega_s$ be the set
of the $s$-th roots of the unity. The term $X_i^r$ takes all the values of $\Omega_s$ equally often.
 \item [P-2]
Let $X^\alpha=X_{j_1}^{\alpha_{j_1}} \cdots X_{j_k}^{\alpha_{j_k}}$ be an interaction term of order $k$ where
$X_{j_i}^{\alpha_{j_i}}$ takes values in $\Omega_{s_{j_i}} $. Let us define $s=\text{lcm}\{s_{j_1}, \ldots, s_{j_k}\}$.
The term $X^\alpha$ takes values in $\Omega_s$ equally often.
\end{enumerate}

Let $X^\alpha$ be a term with  level set $\Omega_s$ on the design $\design$. Let $r_k$ be the number of times $X^\alpha$
takes the value $\omega_k$ on $\fraction$, $k=0,\ldots,s-1$. The polynomial $P(\zeta)$ is associated  to the sequence
$(r_k)_{k=0,\ldots,s-1}$ so that:
\begin{equation*}
P(\zeta) = \sum_{k=0}^{s-1}r_k \zeta^k \qquad \textrm{with } \zeta \in \mathbb C \ .
\end{equation*}
It should be noticed  that
$$E_\fraction(X^\alpha) = \frac1{\# \fraction}\sum_{k=0}^{s-1} r_k \omega_k = \frac1{\# \fraction}P(\omega_1) $$
 See \cite{lang:65} and the Appendix for a review of the properties of such a polynomial $P$.

\begin{proposition}\label{pr:con-orth-Omega} \
Let $X^\alpha$ be a term with level set $\Omega_s$ on  full design $\design$.
\begin{enumerate}
  \item \label{cyclotomic}
 $X^\alpha$  is centered on $\fraction$ if, and only if,
 $$P(\zeta) = \Phi_s(\zeta) \Psi(\zeta)$$ where $\Phi_s$ is the cyclotomic polynomial of the $s$-roots
of the unity and $\Psi$ is a suitable polynomial with integer coefficients.

\item \label{prime} Let $s$ be prime. Therefore, the term $X^\alpha$ is centered on $\fraction$ if, and only if, its $s$
levels appear equally often:
\begin{equation*}
r_0 = \cdots = r_{s-1} = r
\end{equation*}
\item Let $s = p_1^{h_1} \cdots \cdots p_d^{h_d}$, with $p_i$ prime, for $i=1, \ldots, d$. The term $X^\alpha$ is
centered on $\fraction$ if, and only if, the following equivalent conditions are satisfied.
\begin{enumerate}
\item The remainder
\begin{equation*}
H(\zeta)=P(\zeta) \mod \Phi_s(\zeta) \ ,
\end{equation*}
whose coefficients are integer combination of $r_k$, $k= 0, \ldots , s-1$, is identically zero.
\item The polynomial of degree $s$
\begin{equation*}
\tilde P(\zeta)= P(\zeta) \prod_{d | s} \Phi_d(\zeta) \mod (\zeta^s - 1) \ ,
\end{equation*}
whose coefficients are integer combination of the replicates $r_k$, $k= 0, \ldots , s-1$, is identically zero. The indices of the product
are the $d$'s that divide $s$.
\end{enumerate}
\item \label{ort-1} Let $g_i$  be an indicator of a subgroup or of a lateral of a subgroup of $\Omega_s$; i.e.:
$g_i=(g_{i1},\ldots, g_{ij}, \ldots, g_{is})$, $g_{ij}\in \{0,1\}$, such that $\{k \ : \ g_{ik}=1\}$ is a subgroup or a
lateral of a subgroup of $\Omega_s$.

 If the vector of level replicates $(r_0, r_1, \ldots, r_{s-1})$ is a combination with positive weights
of $g_i$:
\begin{equation*}
(r_0, r_1, \ldots, r_{s-1})=\sum a_i \ g_i \quad \textrm{ with } a_i \in \mathbb N
\end{equation*}
$X^\alpha$ is centered.
\end{enumerate}
\end{proposition}
\begin{proof}
\begin{enumerate}
\item \label{item:1} As $\omega_k = \omega_1^k$, the assumption $\sum_k r_k \omega_k = 0$ is equivalent to $P(\omega_1) =
0$. From Section \ref{pr:roots}, Items \ref{it:primit} and \ref{it:cycl-pol}, we know that this implies that
$P(\omega)=0$ for all primitive $s$-roots of the unity, that is, $P(\zeta)$ is divisible by the cyclotomic polynomial
$\Phi_s$.
\item If $s$ is a prime number, the cyclotomic polynomial is $\Phi_s(\zeta) =\sum_{k=0}^{s-1} \zeta^k$. The polynomial
$P(\zeta)$ is divided by the cyclotomic polynomial, and $P(\zeta)$ and $\Phi_s(\zeta)$ have the same degree, therefore
$r_{s-1}
> 0$ and $P(\zeta) = r_{s-1} \Phi(\zeta)$, so that $  r_0 = \cdots = r_{s-1} $.
\item The divisibility shown in Item \ref{item:1} is equivalent to the condition of null remainder. Such a remainder is
easily computed as the reduction of the polynomial $P(\zeta) \mod \Phi_s(\zeta)$. According to the same condition and
Equation (\ref{eq:phiprod}), we obtain that $\tilde P(\zeta)$ is divisible by $\zeta^s - 1$, therefore it also equals $0
\mod \zeta^s-1$.
 \item If $\Omega_p$ is a prime subgroup of $\Omega_s$, then $\sum_{\omega \in \Omega_p} \omega = 0$. Now let us
assume that the replicates  on a primitive subgroup $\Omega_{p_i}$ are 1. Therefore $\sum_{\omega \in \Omega_{p_i}}
\omega=0$ according the equation in Item (\ref{it:s-poly}). The same occurs in the case of the laterals and the sum of
such cases.
 \end{enumerate}
\end{proof}

\medskip \emph{Example} \\
Let us consider the case $s=6$. This situation occurs in the case of mixed-level factorial designs with both three-level
factors and two-level factors. In this case, the cyclotomic polynomial is $\Phi_6(\zeta)=\zeta^2-\zeta+1$ whose roots are
$\omega_1$ and $\omega_5$.  The remainder is
\begin{equation*}
\begin{split}
H(\zeta) & =\sum_{k=0}^{5}r_k \zeta^k  \mod \Phi_6(\zeta) \\& = r_0 + r_1 \zeta + r_2 \zeta^2 + r_3 \zeta^3 + r_4
\zeta^4+ r_5 \zeta^5 \mod\left(\zeta^2-\zeta+1\right) \\
& =  (r_1+r_2-r_4-r_5) \zeta+ (r_0-r_2-r_3+r_5)
\end{split}
\end{equation*}
 The condition $H(\zeta)=0$ implies the following relations concerning the numbers of replicates:
$ r_0+r_1=r_3+r_4 \ , \quad r_1+r_2=r_4+r_5 \ , \quad r_2+r_3=r_0+r_5$, where the first one follows by summing of the
second with the third one. Equivalently:
\begin{equation}\label{rem}
r_0-r_3=r_4-r_1=r_2-r_5 \ .
\end{equation}
Let us consider the replicates corresponding to the sub-group $\{\omega_0, \omega_2, \omega_4\}$ and denote  the
$\min\{r_0,r_2,r_4\}$ by $m_1$. We then consider the replicates corresponding to the lateral of the previous sub-group
$\{\omega_1, \omega_3, \omega_5\}$ and we denote by $m_2$ the $\min\{r_1,r_3,r_5\}$.  We consider the new vector of the
replicates:
\begin{equation*}
\begin{split}
r'&= (r_0',r_1',r_2',r_3',r_4',r_5')\\ &= (r_0-m_1, r_1-m_2, r_2-m_1, r_3-m_2, r_4-m_1, r_5-m_2) \\ &=
r-m_1(1,0,1,0,1,0)-m_2(0,1,0,1,0,1)
\end{split}
\end{equation*}
The vector $r'$ satisfies Equation (\ref{rem}).

As at least  $r_0'$, $r_2'$ or $r_4'$ is zero, the common value in Equation (\ref{rem}) is zero or negative. Moreover, as
at least $r_1'$, $r_3'$ or $r_5'$ is zero,  the common value in Equations (\ref{rem}) is zero or positive. The common
value is therefore zero and $r_0'=r_3'$, $r_1'=r_4'$, $r_2'=r_5'$ and
\begin{equation*}
r' = r_0'(1,0,0,1,0,0)+r_1'(0,1,0,0,1,0)+r_2'(0,0,1,0,0,1)\\
\end{equation*}
A term is therefore centered if the vector of the replicates  is of the form:
\begin{multline*}
(r_0, \ldots, r_5) = a_1 (1,0,0,1,0,0) + a_2 (0,1,0,0,1,0) \\ + a_3 (0,0,1,0,0,1)+ a_4 (1,0,1,0,1,0)+ a_5 (0,1,0,1,0,1)
\end{multline*}
with $a_i$ non negative integers. There are 5 generating integer vectors of the replicate vector.

It should be noticed that if the number of levels of $X^\alpha$ is not prime,  $E_\fraction(X^\alpha)=0$ does not imply
$E_\fraction(X^{r\alpha})=0$. In the previous six-level example, if $X^\alpha$ is centered,  the vector of replicates of
$X^{2\alpha}$ is of the form $(2a_1+a_4+a_5,0,2a_2+a_4+a_5,0,2a_3+a_4+a_5,0)$ and $X^{2\alpha}$ is centered only if
$a_1=a_2=a_3$.

\section{Orthogonal arrays}
\label{sec:project}%
In this sectionn we  discuss the relations between the coefficients $b_\alpha$, $\alpha \in L$, of the counting function
and the property of being an orthogonal array. Let
\begin{equation*}
  \text{OA}(n,s_1^{p_1},\dots,s_m^{p_m},t)
\end{equation*}
be a mixed-level  orthogonal array with $n$ rows and $m$ columns, $m=p_1+\cdots+p_m$, in which $p_1$ columns have $s_1$
symbols, \ldots, $p_k$ columns have $s_m$ symbols, and with strength $t$, as defined e.g. in \cite[p.
260]{wu|hamada:2000}. Strength $t$ means that, for any $t$ columns of the matrix design, all possible combinations of
symbols appear equally often in the matrix.

\begin{definition} \label{de:projectivity}
Let $I$ be a non-empty subset of $\{1,\ldots,m\}$, and let $J$ be its complement set, $J = I^c$.  Let $\design_I$ and $\design_J$ be the
corresponding full factorial designs over the $I$-factors and the $J$-factors, so that $\design = \design_I \times \design_J$.   Let
$\fraction$ be a fraction of $\design$ and let $\fraction_I$ and $\fraction_J$ be its projections.
\begin{enumerate}
\item A fraction $\fraction$ {\em factorially projects} on the $I$-factors if $\fraction_I=s \ \design_I$, that is, the
projection is a full factorial design where each point appears $s$ times.
 \item A fraction $\fraction$ is a {\em mixed orthogonal array} of strength $t$ if it factorially projects on any
 $I$-factors with $\#I=t$.
\end{enumerate}
\end{definition}

Using the notations of Definition \ref{de:projectivity}, for each point $\zeta$ of a complex coded fraction $\fraction$,
we consider the decomposition $\zeta = (\zeta_I, \zeta_J)$ and we denote  the counting function restricted to the
$I$-factors of a fraction by $R_I$, i.e. $ R_I(\zeta_I)$ is the number of points in $\fraction$ whose projection on the
$I$-factors is $\zeta_I$.
\\
We denote  the sub-set of the exponents restricted to the $I$-factors by $L_I$ and  an element of $L_I$ by $\alpha_I$:
\begin{equation*}
L_I= \left\{\alpha_I=(\alpha_1,\ldots,\alpha_j,\ldots,\alpha_m) , \qquad  \alpha_j = 0 \text{ if } j \in J \right\}  \ .
 \end{equation*}
Therefore, for each $\alpha \in L$ and $\zeta \in \design$: $\alpha=\alpha_I+\alpha_J$ and $X^{\alpha}(\zeta)=
X^{\alpha_I}(\zeta_I) X^{\alpha_J}(\zeta_J)$. \\ We denote   the cardinalities of the projected designs by $\# \design_I$
and $\# \design_J$.

\begin{proposition}  \label{pr:projectivity}
\
\begin{enumerate}
\item The number of replicates of the points of a fraction projected onto the $I$-factors is:
\begin{equation*}
R_I(\zeta_I) = \# \design_J \sum_{\alpha_I} b_{\alpha_I} \ X^{\alpha_I}(\zeta_I) \ .
\end{equation*}
\item A fraction \emph{factorially projects onto the $I$-factors} if, and only if,
\begin{equation*}
R_I(\zeta_I)=\# \design_J \ b_0 = \frac {\# \fraction}{\# \design_I} \quad \text{ for all } \ \zeta_I  \ .
\end{equation*}
This is equivalent to all the coefficients of the counting function involving only the $I$-factors being 0:
$$
b_{\alpha_I}=0 \quad \textrm{ with } \ \alpha_I \in L_I, \ \alpha_I \ne (0,0, \ldots ,0) \ .
 $$
 In such a case, the levels of a factor $X_i$, $i \in I$, appear equally often in $\fraction$.
 \item If there exists
a subset $J$ of $\{1,\ldots,m\}$ such that the $J$-factors appear in all the non null elements of the counting function,
 the fraction \emph{factorially projects onto the $I$-factors}, with $I=J^c$.
 \item A fraction is an \emph{orthogonal array of strength $t$}  if, and
 only if, all the coefficients of the counting function up to the order $t$ are zero:
 $$ b_{\alpha}=0 \quad \forall \ \alpha  \textrm{ of order up to }t, \ \alpha \ne (0,0, \ldots ,0)\ .
$$
\end{enumerate}
\end{proposition}

\begin{proof}
\begin{enumerate} \item We obtain:
\begin{equation*}
\begin{split}
 R_I(\zeta_I) &  = \sum_{\zeta_J \in \design_J} R(\zeta_I,\zeta_J) =
 \sum_{\zeta_J \in \design_J}  \sum_{\alpha \in L} b_\alpha \ X^{\alpha}(\zeta_I,\zeta_J) \\
&  =  \sum_{\zeta_J \in \design_J}  \sum_{\alpha \in L} b_\alpha \ X^{\alpha_I}(\zeta_I) X^{\alpha_J}(\zeta_J)\\
&  =  \sum_{\alpha_I \in L_I} b_{\alpha_I} X^{\alpha_I}(\zeta_I) \sum_{\zeta_J \in \design_J} X^{\alpha_J}(\zeta) +
\sum_{\alpha \not \in L_I} b_{\alpha} X^{\alpha_I}(\zeta) \sum_{\zeta_J \in \design_J} X^{\alpha_J}(\zeta)  \ .
\end{split}
\end{equation*}
The thesis follows from $\sum _{\zeta_J \in \design_J} X^{\alpha_J}(\zeta_J) = 0$  if $\alpha_J \ne
(0, 0, \ldots, 0)$ and $\sum _{\zeta_J \in \design_J} X^{\alpha_J}(\zeta_J) = \# \design_J$ if
$\alpha_J =(0, 0, \ldots , 0)$.
 \item The number of replicates of the  points of the fraction
projected onto the $I$-factors, $ R_I(\zeta_I) = \# \design_J \sum_{\alpha_I} b_{\alpha_I} \ X^{\alpha_I}(\zeta_I)$, is a
polynomial and it is a constant if all the coefficients $b_{\alpha_I}$, with $\alpha_I \neq (0, 0, \ldots , 0)$, are
zero.
 \item This condition implies that the $b_{\alpha_I}$'s are zero, if $\alpha_I \neq (0,
0, \ldots , 0)$, and the thesis follows from the previous item.
 \item This item follows from the previous items and the definition.
\end{enumerate}
\end{proof}

\medskip\emph{Remarks}
\begin{enumerate}
 \item If a fraction factorially projects onto the $I$-factors,  its cardinality  must be equal to, or a multiple of the
cardinality of $\design_I$.
 \item If the number of levels of each factors is a prime, the condition $b_{\alpha_i}=0$ for each $i\in I$ and $0<\alpha_i\leq n_i-1$
in  Items (2) and (3) of the previous Proposition, simplify to $E_\fraction\left(X_i\right)=0$, according to Item
(\ref{prime}) of  Proposition \ref{pr:con-orth-Omega}.

 \end{enumerate}

\section{Regular fractions: a partial generalization to mixed-level design} \label{sec:regular}

A short  review of the theory of regular fractions is here made from the view point of the present paper. Various
definitions of regular fraction appear in literature, e.g. in the books by \cite[p. 123]{raktoe|hedayat|federer:81},
\cite[p. 125]{collombier:96}, \cite[p. 70]{kobilinsky:97}, \cite[p. 164]{dey|mukerjee:1999}, \cite[p.
305]{wu|hamada:2000}. To our knowledge, all the definitions are known to be equivalent if all the factors have the same
prime number of levels, $n=p$. The definition based on Galois Field computations is given for $n=p^s$ power of a prime
number. All definitions assume symmetric factorial design, i.e. all the factors have the same number of levels.

Regular fraction designs are usually considered for qualitative factors, where the  coding of the levels is arbitrary.
The integer coding, the GF$(p^s)$ coding, and the roots of the unity coding, as introduced by \cite{bailey:1982} and
used extensively in this paper, can all be used. Each of those codings is associated to specific ways of characterizing
a fraction, and even more important for us, to a specific basis for the responses. One of the possible definitions of a
regular fraction refers to the property of non-existence of partial confounding of simple and interaction terms, and
this property has to be associated to a specific basis, as explicitly  pointed out in \cite{wu|hamada:2000}.

In our approach, we use polynomial algebra with complex coefficients, the $n$-roots of the unit coding, and the idea of
indicator polynomial function, and we make no assumption about the number of levels. In the specific coding we use, the
indicator polynomial is actually a linear combination of monomial terms which are centered and orthogonal on the full
factorial design. We refer to such a basis to state the no-partial confounding property.

The definition of the regular fraction is hereafter generalized in the symmetric case with a prime number of levels. The
new setting includes asymmetric design with any number of levels.    Proposition \ref{th:reg} below does not include
regular fractions defined in GF$(p^s)$. A full discussion of this point shall be published elsewhere.

We consider a fraction \emph{without replicates}. Let $n=\text{lcm}\{n_1,\ldots,n_m\}$. It should be recalled that
$\Omega_n$ is the set of the $n$-th roots of the unity, $\Omega_n = \{\omega_0,\ldots,\omega_{n-1}\}$.  Let $\mathcal L$
be a subset of exponents, $\mathcal L \subset L = \mathbb Z_1\times\cdots\times\mathbb Z_m$, containing $(0,\ldots,0)$
and let $l$ be its cardinality ($l>0$). Let $e$ be a map from $\mathcal L$ to $\Omega_n$, $e : \mathcal L \rightarrow
\Omega_n$.

\begin{definition} \label{def:reg}
A fraction $\fraction$ is \emph{regular} if
\begin{enumerate}
 \item $\mathcal L$ is a sub-group of $L$,
 \item $e$ is a group homomorphism, $e([\alpha+\beta])= e(\alpha ) \ e(\beta)$ for each $\alpha, \beta \in \mathcal L$,
 \item the equations
\begin{equation} \label{eq:def-eq}
X^\alpha = e(\alpha) \quad , \qquad  \alpha \in \mathcal L
\end{equation}
define the fraction $\fraction$, i.e. they are a set of generating equations, according to Section \ref{sec:gen-eq}.
Equations (\ref{eq:def-eq}) are also called the \emph{defining} equations of $\fraction$.
\end{enumerate}
If $\mathcal H$ is a minimal generator of the group $\mathcal L$,  Equations $X^\alpha = e(\alpha)$,  $\alpha \in
\mathcal H \subset \mathcal L$, are called a minimal set of \emph{generating} equations.
\end{definition}
 It should be noticed that our situation is general because the values $e(\alpha)$ can be  different from $1$.
From items (1) and (2) it follows that a necessary condition is that the $e(\alpha)$'s must belong to the subgroup
spanned by the values of $X^\alpha$. For example, for $n_1=n_2=n=6$, an equation such as $X_1^3 X_2^3 = \omega_2$ cannot
be a defining equation.

For example, in the fraction of Section \ref{sec:gen-eq}, we have: $\mathcal H=\left\{(1,1,2,0), (1,2,0,1)\right\}$ and
$e(1,1,2,0)=e(1,2,0,1)=\omega_0=1$. The set $\mathcal L$ is:  $\{(0,0,0,0),(0,1,1,2), $ \\ $(0,2,2,1),(1,1,2,0),
(2,2,1,0), (1,2,0,1), (2,1,0,2), (1,0,1,1), (2,0,2,2)\}$.

\begin{proposition} \label{th:reg}
Let $\fraction$ be a fraction.  The following statements are equivalent:
\begin{enumerate}
 \item \label{it:reg} Fraction $\fraction$ is regular according to Definition \ref{def:reg}.
  \item \label{it:reg-F} The indicator function of the fraction has the form
\begin{equation*}
F(\zeta)=\frac 1 l \sum_{\alpha \in \mathcal L} \overline{e(\alpha)}\ X^\alpha(\zeta) \qquad \zeta \in \design
\end{equation*}
where $\mathcal L$ is a given subset of $L$ and $e: \mathcal L \to \Omega_n$ is a given mapping.
 \item \label{it:reg-con} For each $\alpha, \beta \in  L$, the parametric functions represented on $\fraction$ by the terms $X^\alpha$
 and $X^\beta$ are either orthogonal or totally confounded.
\end{enumerate}
\end{proposition}
\begin{proof}

First we prove the equivalence between (\ref{it:reg}) and (\ref{it:reg-F}).

(\ref{it:reg}) $\Rightarrow$ (\ref{it:reg-F}).
\\
Let $\fraction$ be a regular fraction and let $X^\alpha=e(\alpha)$ be its defining equations with $\alpha \in \mathcal
L$, $\mathcal L$ a sub-group of $L$ and $e$ a homomorphism.

If, and only if, $\zeta \in \fraction$:
\begin{equation*}
\begin{split}
0  &=
 \sum_{\alpha \in \mathcal L} \vert X^\alpha(\zeta)-e(\alpha) \vert^2 =
 \sum_{\alpha \in \mathcal L} \left( X^\alpha(\zeta)-e(\alpha) \right) \overline {\left( X^\alpha(\zeta)-e(\alpha) \right)} \\
 & = \sum_{\alpha \in \mathcal L} \left( X^\alpha(\zeta) \overline  {X^\alpha(\zeta)} + e(\alpha) \overline  {e(\alpha)}
 - e(\alpha) \overline  {X^\alpha(\zeta)} - \overline  {e(\alpha)} X^\alpha(\zeta)\right)\\
 & = 2 \ l -\sum_{\alpha \in \mathcal L} \overline{\overline{e(\alpha)}  {X^{\alpha}(\zeta)}}- \sum_{\alpha \in \mathcal L} \overline{e(\alpha)}
  \ X^\alpha(\zeta) = 2\left( l -  \sum_{\alpha \in \mathcal L} \overline{e(\alpha)} \ X^\alpha(\zeta)\right)
\end{split}
\end{equation*}
therefore
\begin{displaymath}
\frac 1 l \sum_{\alpha \in {\mathcal L}} \overline{e(\alpha)} \ X^\alpha(\zeta) - 1=0 \quad \textrm{if, and only if, }
\quad \zeta \in \fraction \ .
\end{displaymath}

The function $F =\frac 1 l \sum_{\alpha \in {\mathcal L}} \overline{e(\alpha)} \ X^\alpha$ is an indicator function, as
it can be shown that $F=F^2$ on $\design$. In fact, $\mathcal L$ is a sub-group of $L$ and $e$ is a homomorphism;
therefore:
\begin{equation*}
\begin{split}
F^2 &=\frac 1 {l^2} \sum_{\alpha \in {\mathcal L}} \sum_{\beta \in {\mathcal L}} \overline{e(\alpha) \ e(\beta)} \
X^{[\alpha+\beta]}= \frac 1 {l^2} \sum_{\alpha \in {\mathcal L}}
\sum_{\beta \in {\mathcal L}} \overline{e([\alpha+\beta])} \ X^{[\alpha+\beta]}= \\
  &=
\frac 1 {l^2} \sum_{\gamma \in {\mathcal L}} l \ \overline{e(\gamma)} \ X^\gamma =F
\end{split} \ .
\end{equation*}
It follows that $F$ is the indicator function of $\fraction$, and $ b_\alpha =\frac {\overline{e(\alpha)}} l$, for all
$\alpha \in {\mathcal L}$.

(\ref{it:reg-F}) $\Rightarrow$ (\ref{it:reg}).\\
It should be noticed that an indicator function is real valued, therefore $\overline F = F$.
\begin{equation*}
\frac 1 l \sum_{\alpha \in \mathcal L} \vert X^\alpha(\zeta) - e(\alpha) \vert^2 = 2 - \overline{ F(\zeta)}-F(\zeta)= 2 - 2 F(\zeta) =
\begin{cases}
0 & \text{on $\fraction$} \\
2  & \text{on $\design \setminus \fraction$} \ .
\end{cases}
\end{equation*}
Equations  $X^\alpha = e(\alpha)$, with $\alpha \in \mathcal L$, define the fraction $\fraction$ as the generating
equations of a regular fraction. It is easy to see that $\mathcal L$ is a group. In fact, if $\gamma=[\alpha+\beta]\notin
\mathcal L$, there exists one $\zeta$ such that $X^{\gamma}(\zeta)=X^{\alpha}(\zeta)X^{\beta}(\zeta)=e(\alpha)e(\beta)
\subset \Omega_n$ and the value $e(\alpha)e(\beta)$  only depends on $\gamma$. By repeating the previous proof,  the
uniqueness of the polynomial representation of the indicator function leads a contradiction.

Now we prove the equivalence between (\ref{it:reg-F}) and (\ref{it:reg-con}).

(\ref{it:reg-F}) $\Rightarrow$ (\ref{it:reg-con})\\ The non-zero coefficients of the indicator function are of the form
$b_\alpha=\overline {e(\alpha)} / l$.

We consider two  terms $X^\alpha$ and $X^\beta$ with $\alpha, \beta \in L$. If $[\alpha-\beta] \notin \mathcal L$ then
$X^\alpha$ and $X^\beta$ are orthogonal on $\fraction$ as the coefficient $b_{[\alpha-\beta]}$ of the indicator function
equals 0. If $[\alpha-\beta] \in \mathcal L$ then $X^\alpha$ and $X^\beta$ are confounded because
$X^{[\alpha-\beta]}=e([\alpha-\beta])$; therefore $X^\alpha=e([\alpha-\beta])\ X^\beta$.

(\ref{it:reg-con}) $\Rightarrow$ (\ref{it:reg-F}).
\\
 Let $\mathcal L$ be the set of exponents of the  terms confounded with a constant:
\begin{equation*}
\mathcal L = \left\{ \alpha \in L : X^\alpha = \text{constant} = e(\alpha), \quad e(\alpha) \in \Omega_n \right\} \ .
\end{equation*}
For each $\alpha \in \mathcal L$, $b_\alpha = \overline{e(\alpha)}\ b_0$. For each $\alpha \notin \mathcal L$, because of
the assumption, $X^\alpha$ is orthogonal to $X^0$, therefore $b_\alpha=0$.
\end{proof}

\begin{corollary}
Let $\fraction$ be a regular fraction with $X^\alpha=1$ for all the defining equations.   $\fraction$ is therefore
self-conjugate and a multiplicative subgroup of $\design$.
\end{corollary}

\begin{proof}
It follows from Prop. \ref{pr:ort-bc-alpha} Item \ref{ort-4}.
\end{proof}

The following proposition extends a result presented in \cite{fontana|pistone|rogantin:2000} for the two level case.
\begin{proposition} \label{pr:reg-cont}
Let $\fraction$ be a fraction with indicator function $F$. We denote the set of the exponents $\alpha$  such that $\frac
{b_\alpha}{b_0} = \overline{e(\alpha)} \in \Omega_n$ by $\mathcal L$ . The indicator function can be written as
\begin{equation*}
F(\zeta) = b_0 \sum_{\alpha \in \mathcal L} \overline{e(\alpha)}\ X^\alpha(\zeta) + \sum_{\beta \in \mathcal K} b_\beta \ X^\beta(\zeta)
\qquad \zeta \in \fraction \ , \ \mathcal L \cap \mathcal K = \emptyset \ .
\end{equation*}
It follows that $\mathcal L$ is a subgroup and the equations $X^\alpha =  e(\alpha)$, with $\alpha \in \mathcal L$, are
the defining equations of the smallest regular fraction $\fraction_{r}$  containing $\fraction$ restricted to the factors
involved in the $\mathcal L$-exponents.
\end{proposition}

\begin{proof}
The coefficients $b_\alpha$, $\alpha \in \mathcal L$, of the indicator function $F$ are of the form  $b_0
\overline{e(\alpha)}$. Therefore, from the extremality of $n$-th roots of the unity,  $X^\alpha(\zeta)=e(\alpha)$ if
$\zeta \in \fraction$ and $X^\alpha(\zeta)F(\zeta)=e(\alpha)F(\zeta)$ for each $\zeta \in \design$ and $\mathcal L$ is a
group.

We denote  the indicator function of $\fraction_r$ by $F_r$. For each $\zeta \in \design$ we have:
\begin{equation*}
\begin{split}
F(\zeta) F_r(\zeta)&= \frac 1 l F(\zeta) \sum_{\alpha \in \mathcal L} \overline{e(\alpha)}\ X^\alpha(\zeta) =
\frac 1 l  \sum_{\alpha \in \mathcal L} \overline{e(\alpha)}\ X^\alpha(\zeta) \ F(\zeta)= \\
&= \frac 1 l  \sum_{\alpha \in \mathcal L} \overline{e(\alpha)}\ e(\alpha) \ F(\zeta) = \frac 1 l \ l \ F(\zeta) \ .
\end{split}
\end{equation*}
The relation  $F(\zeta) F_r(\zeta)=F(\zeta)$ implies $\fraction \subseteq \fraction_r$. The fraction $\fraction_r$ is
minimal because we have collected all the terms confounded with a constant.
\end{proof}

\subsection*{Remark}
Given generating equations $X^{\alpha_1}=1, \ldots, X^{\alpha_h}=1$, with $\{\alpha_1, \ldots, \alpha_h\}=\mathcal H
\subset \mathbb Z_{n_1} \times \cdots \times \mathbb Z_{n_m}$, the corresponding fraction $\fraction$ is a subgroup of
$\Omega_{n_1} \times \cdots \times \Omega_{n_m}$. If the same fraction is represented in the additive notation, such a
set of treatment combinations is the principal block of a single replicate generalized cyclic design, see
\cite{john|dean:75}, \cite{dean|john:75} and \cite{lewis:79}. A complex vector of the form
\begin{equation*}
\left( e^{\jei 2\pi \frac{k_1}{n_1}},e^{\jei 2\pi \frac{k_2}{n_2}},\ldots,e^{\jei 2\pi \frac{k_m}{n_m}}\right) \qquad
\textrm{with } \ 0 \leq k_i < n_i, \ i=1,\ldots,m \ ,
\end{equation*}
is in fact a solution of the generating equations if, and only if,
\begin{equation} \label{eq:gen-treat}
\left\{
\begin{array}{ll}
\sum_{j=1}^m \alpha_{ij} \ \gamma_j \ k_j=0  \mod s &  \\
0 \leq k_j < n_j
\end{array} \right.
\end{equation}
with $s=\textrm{lcm}\{n_1,\ldots,n_m\}$ and $\gamma_j=\frac s{n_j}$.

A set of generators can be computed from  Equation (\ref{eq:gen-treat}). It should be noticed that the following
equivalent integer linear programming problem does not involve computation mod $s$, see \cite{schrijver:86}
\begin{equation} \label{eq:gen-treat2}
\left\{
\begin{array}{ll}
\sum_{j=1}^m \alpha_{ij} \ \gamma_j \ k_j-sq=0  \\
0 \leq k_j < n_j \ , \ q \ge 0
\end{array} \right.
\end{equation}
In \cite{lewis:82}, the monomial part of our defining equations is called defining contrast, according to
\cite{bailey|gilchrist|patterson:77}. The paper contains extensive tables of the generator subgroups of the treatment
combinations and the corresponding defining contrasts.

Viceversa, given a set of generators  of the treatment combinations,
 \begin{equation*}
\left\{b_1,\ldots,b_r \ | \ b_i=\left(b_{i1},\ldots ,b_{im}\right) \right\} \ ,
 \end{equation*}
  Equation (\ref{eq:gen-treat}) with indeterminates $\alpha_i$
\begin{equation*}
\left\{
\begin{array}{ll}
\sum_{j=1}^m  \gamma_j \ b_{ij}  \ \alpha_{j} =0 \mod s &  \\
0 \leq \alpha_{j} < n_j
\end{array} \right.
\end{equation*}
produces generating equations for the fraction.

\section{Examples} \label{sec:ex}

\medskip \textbf{A regular fractions with $\mathbf{n=3}$.}\\ Let us consider the classical $3^{4-2}$ fraction of Section
\ref{sec:gen-eq}. Its indicator function is:
\begin{eqnarray*}
F & = & \frac 1 9 \left(
1 + X_2 X_3 X_4 + X_2^2 X_3^2 X_4^2 + X_1 X_2 X_3^2 + X_1^2 X_2^2 X_3   \right. \\
& &  \left. + X_1 X_2^2 X_4 + X_1^2 X_2 X_4^2 + X_1 X_3 X_4^2 + X_1^2 X_3^2 X_4 \right) \ .
\end{eqnarray*}
We can observe that the coefficients are all equal to $\frac 1 9 $.  The minimum order of interactions that appear in the
indicator function is 3, therefore the fraction is an orthogonal array of strength 2. All the defining equations are of
the form $X^\alpha=1$, therefore the fraction is self-conjugate.

\medskip \textbf{A regular fraction with $\mathbf{n=6}$.}\\ Let us consider a $6^3$ design. From property [P-2] of Section
\ref{sec:orth}, the terms $X^\alpha$ take values either in $\Omega_{6}$ or in one of the two subgroups either $\{1,
\omega_3\}$ or $\{1, \omega_2, \omega_4 \}$.\\ Let $\fraction$ be a fraction whose generating equations are: $ X_1^3
X_2^3 X_3^3=\omega_3$ and $ X_2^4X_2^4X_3^2=\omega_2$. In this case  we have: $\mathcal H=\left\{(3,3,3),
(4,4,2)\right\}$ and $e(3,3,3)=\omega_3$, $e(4,4,2)=\omega_2$. The set $\mathcal L$ is:
 $
\left\{  (0,0,0),(3,3,3),(4,4,2), (2,2,4), (1,1,5), (5,1,1) \right\}$.  The full factorial design has 216 points and the
fraction has 36 points. The indicator function is:
\begin{equation*}
F  =  \frac 1 6 \left( 1 +\omega_3 X_1^3 X_2^3 X_3^3 +\omega_4 X_1^4 X_2^4 X_3^2 + \omega_2 X_1^2 X_2^2 X_3^4 + \omega_1
X_1 X_2 X_3^5 + \omega_5 X_1^5 X_2 X_3 \right)
\end{equation*}
It should be noticed that this fraction is an $OA(36,6^3,2)$.

\medskip \textbf{An $\mathbf{OA(18,2^1 3^7,2)}$.}\\ We consider the  fraction  of a $2 \times 3^7$ design with 18 runs,  taken from
\cite[Table 7C.2]{wu|hamada:2000} and recoded with complex levels. Here $X_1$ takes values in $\Omega_2$, $X_i$, with
$i=2, \ldots, 8$, and their interactions take values in $\Omega_3$, and the interactions involving $X_1$ take values in
$\Omega_6$.
\\
All the 4374 $X^\alpha$ terms  of the fraction have been computed in SAS using $\mathbb Z_2$, $\mathbb Z_3$ and $\mathbb
Z_6$ arithmetic.  The replicates of the values in the relevant $\mathbb Z_k$ have then been computed for each terms. We
found:
\begin{enumerate}
\item 3303 centered responses. These are characterized by Proposition \ref{pr:con-orth-Omega}.   The replicates are of
the type: $(9,9)$, $(6,6,6)$, $(3,3,3,3,3,3)$ and $(9,0,0,9,0,0)$. We have:
\begin{enumerate}
 \item the two-level simple term and 1728 terms involving only the three-level factors (14 of order 1,  84 of order 2, 198 of order 3,
422 of order 4, 564 of order 5, 342 of order 6 and 104 of order 7);
 \item 1574 terms involving both the two-level factor and the three-level factors (14 of
order 2, 66 of order 3,  188 of order 4, 398 of order 5, 492 of order 6, 324 of order 7 and 92 of order 8).
\end {enumerate}
\item 9 terms with corresponding $b_\alpha$ coefficients equal to $b_0=\frac {18}{2 \times 3^7}=3^{-5}$; \item 1062 terms
with corresponding coefficients different from zero and $b_0$: 450 terms involving only the three-level factors (80 of
order 3, 138 of order 4, 108 of order 5, 100 of order 6 and 24 of order 7) and 612 terms involving both the two-level
factor and the three-level factors (18 of order 3,  92 of order 4, 162 of order 5, 180 of order 6, 124 of order 7 and 36
of order 8).
\end {enumerate}
Some statistical properties of the fraction are:
\begin{enumerate}
\item Analyzing the centered responses we can observe that:
\begin{enumerate}
 \item All the 15 simple terms are centered.
\\
All the 98 interactions of order 2 (84 involving only the three-level factors and 14 also involving  the two-level
factor) are centered.  This implies that both the ``linear''  terms and the ``quadratic'' terms of the three-level
factors are mutually orthogonal and they are orthogonal to the two-level factor.
\\
The fraction is a mixed orthogonal array of strength 2.
 \item The fraction factorially projects onto the following factor subsets:
\begin{eqnarray*}
& & \{X_1, X_2, X_3\}, \ \{X_1, X_2, X_4\}, \ \{X_1, X_2, X_5\}, \ \{X_1, X_2, X_6\}, \\
& & \{X_1, X_3, X_6\}, \ \{X_1, X_3, X_7\}, \ \{X_1, X_4, X_5\}, \ \{X_1, X_4, X_8\}, \\
& & \{X_1, X_5, X_8\}, \ \{X_1, X_6, X_7\}, \ \{X_1, X_6, X_8\} \ .
\end{eqnarray*}
All the terms of order 1, 2 and 3 involving the same set  of factors are in fact centered. \item The minimal regular
fraction containing our fraction restricted to the three-level factors has the following defining relations:
\begin{eqnarray*}
& & X_2^2 X_4^2 X_5=1 \ , \  X_2 X_4 X_5^2=1 \ , \\ & &
 X_2 X_3 X_4^2 X_6 X_7 X_8=1 \ , \  X_2^2 X_3^2 X_4 X_6^2 X_7^2 X_8^2=1 \ , \\ & &
  X_2^2 X_3 X_5^2 X_6 X_7 X_8=1 \ , \  X_2 X_3^2 X_5 X_6^2 X_7^2 X_8^2=1 \ , \\ & &
  X_3 X_4 X_5 X_6 X_7 X_8=1 \ , \    X_3^2 X_4^2 X_5^2 X_6^2 X_7^2 X_8^2=1  \ .
\end{eqnarray*}
\item The non centered terms have  levels in $\Omega_6$ and in $\Omega_3$.
\end {enumerate}
\end {enumerate}

\subsection*{Acknowledgments}
 We wish to thank  many colleagues for their helpful and interesting comments, especially G.-F. Casnati, R.
Notari, L. Robbiano, E. Riccomagno, H.P. Wynn and K.Q. Ye. Last but not least, we extensively  used the comments and
suggestions made by the anonymous referees of the previous versions. We regret we are unable to thank them by name.

\section{Appendix: Algebra of the $n$-th roots of the unity.}  \label{pr:roots} %
We hereafter list  some facts concerning the algebra of the complex $n$-th roots of the unity,  for ease of reference.
\begin{enumerate}
\item  \label{it:conj} The conjugate of a $n$-th root of the unity equals its inverse:
 $\overline{\omega_k}= \omega_k^{-1} = \omega_{[-k]}$ for all $\omega_k \in \Omega_n$.
 \item \label{it:pow} If $\zeta \ne \omega_m$, we obtain: $ \prod_{k=0 \ k\ne m}^{n-1}  (\zeta - \omega_k)=\frac{\zeta^n-1}{\zeta-\omega_m}= \sum_{h=0 }^{n-1} \omega_m^{n-h-1} \
\zeta^h $ where the last equality follows from  algebraic computation. Therefore, for $\zeta=\omega_m$:
\begin{equation*}\prod_{k=0 \ k\ne m}^{n-1}  (\omega_m - \omega_k) =  \sum_{h=0 }^{n-1} \omega_m^{n-h-1} \ \omega_m^h = n \ {\omega_m}^{n-1}=n \ \overline{\omega_m}
\end{equation*}
and especially: $\prod_{k=1 }^{n-1} (1- \omega_k) =n $.

\item \label{it:s-poly} We have: $ \zeta^n-1 =  (\zeta-\omega_0) \cdots (\zeta-\omega_{n-1}) = \sum_{k=0}^{n-1}
(-1)^{n-k} S_{n-k}\left(\omega_0,\ldots,\omega_{n-1}\right) \zeta^k $ where $S_{n-k}\left(x_0,\ldots,x_{n-1} \right)$ is
the elementary symmetric polynomial  of order $n-k$. We therefore obtain the following notable cases:
\begin{itemize}
\item[-] $S_{1}\left(\omega_0,\ldots,\omega_{n-1}\right)=\sum_k \omega_k =  0$
 \item[-]
$S_{2}\left(\omega_0,\ldots,\omega_{n-1}\right)=\sum_{{\ell}<m}  \omega_{\ell} \ \omega_m=0$ \item[-]
$S_{n}\left(\omega_0,\ldots,\omega_{n-1}\right)=\prod_k \omega_k = (- 1)^{n+1}$
\end{itemize}
where the indices of the sums and products are from 0 to $n-1$.
 \item \label{it:primit} Let $\omega$ be
a primitive $n$-th root of the unity, that is, a generator of $\Omega_n$ as a cyclic group: $ \left\{1, \omega, \omega^2,
\ldots, \omega^{n-1}\right\} = \Omega_n   \ . $
\\
The root $\omega_p \in \Omega_n$ is primitive if $p$ is relatively prime with $n$. In particular, $\omega_1$ is a
primitive root and, for $\omega_k \in \Omega_n$, we obtain: $\omega_k=(\omega_1)^k$. If $n$ is a prime number, all the
roots of the unity, except 1, are primitive roots. The number of the primitive $n$-th roots of the unity is denoted by
$\phi(n)$.

\item \label{it:cycl-pol} Given an algebraic number $x$, the  unique irreducible monic polynomial  of the smallest degree
with rational coefficients $P$ such that $P(x)=0$ and whose leading coefficient is 1, is called  the minimal polynomial
of  $x$. The minimal polynomial of a primitive $n$-th root of the unity  is called the cyclotomic polynomial
$\Phi_n(\zeta)$ and its degree is $\phi(n)$:
\begin{equation*} \Phi_n(\zeta) = \prod_{p} (\zeta - \omega_p) \ , \quad \zeta \in \mathbb C,
\ \omega_p \in \Omega_n \text{ primitive } n \text{-th root of the unity.}
\end{equation*}
If $n$ is prime, the minimal polynomial of  a primitive $n$-th root of the unity is $ \Phi_n(\zeta)
=\zeta^{n-1}+\zeta^{n-2} + \cdots + 1 $. Moreover:
\begin{equation} \label{eq:phiprod}
\zeta^n-1= \Phi_n(\zeta) \cdot \cdots \Phi_d(\zeta) \cdot \cdots \Phi_1(\zeta)\qquad \text{where } d \text{ divides }n
 \ .
\end{equation}
\item
 The recoding in Equation (\ref{map}) is a polynomial  function of degree $n-1$ and complex coefficients in both directions:
\begin{eqnarray} \label{recoding}
\omega_k &=& \sum_{s=0}^{n-1} \ \omega_s \   \frac{\prod_{h=0,h\ne s}^{n-1}\left(x -  h\right)}{\prod_{h=0,h\ne s}^{n-1}\left(s - h\right)} \ , \qquad  x=k \in \{0,\ldots,n-1\}  \nonumber   \\
k& =& \frac 1 n \sum_{h=0}^{n-1} \zeta^h  \ \sum_{s=1}^{n-1}  s \ \omega_{[s-sh]} \ \quad \ , \qquad \zeta=\omega_k \in \Omega_n \ .
\end{eqnarray}
The last  Equation  follows from
\begin{equation*}
k=  \sum_{s=1}^{n-1} s \  \frac{\prod_{h=0,h\ne s}^{n-1}\left(\zeta -  \omega_h\right)}{\prod_{h=0,h\ne s}^{n-1}\left(\omega_s -
\omega_h\right)}  \quad \ , \qquad \zeta=\omega_k \in \Omega_n
\end{equation*}
and from the properties of the $n$-th roots of the unity, see  Item \ref{it:pow}.
\end{enumerate}

\end{document}